\documentclass{amsart}
\usepackage{latexsym,amsxtra,amscd,ifthen}
\usepackage{amsfonts}
\usepackage{verbatim}
\usepackage{amsmath}
\usepackage{amsthm}
\usepackage{amssymb}

\numberwithin{equation}{section}

\theoremstyle{plain}
\newtheorem{theorem}{Theorem}[section]
\newtheorem{theorema}[theorem]{Theorem}

\newtheorem{prop}[theorem]{Proposition}
\newtheorem{lemma}[theorem]{Lemma}
\newtheorem{cor}[theorem]{Corollary}

\newtheorem{rema}[theorem]{Remark}

\newtheorem{quest}[theorem]{Question}

\theoremstyle{definition}
\newtheorem{defi}[theorem]{Definition}

\newcommand{\N}{\ensuremath{\mathbb N}}

\def\la{\langle}
\def\ra{\rangle}

\def\deg{\mbox{deg\,}}
\def\dim{\mbox{dim\,}}

\def\gkd{\mbox{GK--dim\,}}

\def\qcoh{\mbox{q-coh\,}}

\def\tor{\mbox{Tor\,}}

\def\Zpl{\mbox{\bf Z}_+}

\begin{document}

\title{Sets of Hilbert series and their applications}
\author{Dmitri Piontkovski}

      \address{ Central Institute of Economics and Mathematics\\
                      Nakhimovsky prosp. 47, Moscow 117418,  Russia}
\thanks{Partially
supported by the grant 02-01-00468 of the Russian Basic
Research Foundation}

\email{piont@mccme.ru}

\begin{abstract}
We consider graded finitely presented algebras and modules
over a field. Under some restrictions,  the set of Hilbert
series of such algebras (or modules) becomes finite. Claims
of that types imply rationality of Hilbert and Poincare
series of some algebras and modules, including periodicity
of Hilbert functions of common (e.g., Noetherian) modules
and algebras of linear growth.
\end{abstract}



\maketitle

\section{Introduction}

We consider graded finitely presented algebras and modules
over a fixed basic field $k$. The set of Hilbert series of
such algebras (or modules), that satisfy some additional
restrictions, becomes finite. We give several applications
of the claims of that type: the rational dependence of
Hilbert series and rationality of Poincare series of
ideals in finitely presented algebras, and
 the periodicity of
Hilbert functions of common finitely presented algebras and
modules of linear growth.


The paper
is organized as follows. In subsection~\ref{ss_notations}
we introduce some notations. Then, in
section~\ref{s_sets_Hilbert}, we present our key results
about sets of Hilbert series. In particular, we prove here
the following

\begin{theorem}[Theorem~\ref{main}, $(a)$]
\label{t-main-intro}
 Given 4 positive integers $n,a,b,c$, let
$D(n,a,b,c)$ denote the set of all connected graded
algebras $A$ over a fixed field $k$ with at most $n$
generators such that $m_1 (A) \le a, m_2(A) \le b$, and
$m_3(A) \le c$. Then the set of Hilbert series of algebras
from $D(n,a,b,c)$ is finite.
\end{theorem}

Here $m_i (A) = \sup \{ j | \tor_i^A(k, k)_j \ne 0 \}$ (if
$\tor_i^A(k, k) =0$, we put $m_i (A) =0$); in particular,
$m_1(A)$ is the exact bound for degrees of generators of
the algebra $A$, and $m_2(A)$ is the exact bound for
degrees of relations of $A$. For example, an algebra $A$ is
Koszul iff
 $m_i (A) \le i$ for all $i \ge 0$.

This Theorem~\ref{t-main-intro} has been also proved
in~\cite{prepr} (a version for Koszul algebras had been
early proved in~\cite{pp}), but here we give another proof.
This new proof seems more clear and elementary. It is based
on the following

\begin{theorema}[Theorem~\ref{an-mod}]
\label{t-an-mod-intro}
Let $D \ge D_0$ be two integers. Consider a set $E(D_0, D)$
of all graded (bi)modules $M$ over  connected graded
algebras $A$ such that  $A$ and $M$ are generated by $\le
D$ elements, the generators of $M$ have degrees at least
$D_0$, and the generators and relations of both $A$ and $M$
have degrees at most $D$. Then the set of all Hilbert
series of (bi)modules from $E(D_0, D)$ has no infinite
ascending chains (with respect to the lexicographical order
on Hilbert series).
\end{theorema}

It is a generalization of a result of D. Anick~\cite{an4},
where $M=A$.

The above theorems gives the following new

\begin{cor}[Corollary~\ref{mod}]
\label{c-modules_best-intro} \label{mod-intro}
Let $D > 0$ be an integer, and let $A$ be a finitely presented
algebra. Then the set of Hilbert series of right-sided ideals in
$A$ having generators and relations in degrees at most $D$ is
finite.
%
\end{cor}

As before, here we put $m_i (M) = \sup \{ j | \tor_i^A(M,
k)_j \ne 0 \}$ (with $m_i (M) =0$ if $\tor_i^A(M, k) =0$);
in particular, $m_0(M)$ is the exact bound for degrees of
generators of  $M$, and $m_1(M)$ is the exact bound for
degrees of relations of $M$.

Our results on sets of Hilbert series gives some
interesting applications~\cite{prepr}. A finitely presented
graded module $M$ over a connected graded algebra $A$ is
called {\it effectively coherent} if there is a function
$D_M : \N \to \N$ such that, whenever a graded submodule $L
\subset M$ is generated in degrees $\le d$, the relations
of $L$ are concentrated in degrees at most $D(d)$. A module
$M$ is called {\it effective for series} if for every
integer $d$ there is only a finite number of possibilities
for Hilbert series of submodules of $M$ generated in
degrees at most $d$. Theorem~\ref{t-main-intro} has been
essentially used to establish the following

\begin{theorem}[\cite{prepr}]
\label{t-intro-str-noet}
(a) Every strongly Noetherian connected algebra over an
algebraically closed field is effectively coherent.

(b) Every effectively coherent algebra is effective for
series.
\end{theorem}

Recall that an algebra $A$ is called strongly
Noetherian~\cite{asz} if an algebra $A\otimes C$ is
Noetherian for every Noetherian commutative $k$--algebra
$C$; in particular, the most of common rings of
non-commutative projective geometry are strongly
Noetherian~\cite{asz}.

Here we consider other applications of sets of Hilbert
series. First, we consider (in section~\ref{s-periodic})
 the following question:
when  Hilbert function $f_V(n):= \dim V_n$ of a graded
algebra or a module $V$ is periodic? An obvious necessary
condition is that the algebra (module) $V$ must have linear
growth, that is, $\gkd V \le 1$. It happens that in some
common cases  this condition is sufficient.

\begin{theorem}[Theorem~\ref{t-main-GK1}]
\label{t-main-GK1-intro}
 Let $M$ be a finitely presented graded module
    over a connected finitely presented algebra $A$.
Suppose that $\gkd M =1$ and  at least one of the following
conditions holds:

  (a)  the field $k$ is finite;

  (b) the vector spaces
  $\tor_2^A (M,k)$  and $\tor_3^A (k,k)$ have finite dimensions.

    Then the Hilbert series $M(z)$ is rational, that is, the Hilbert function $f_M(n) = \dim M_n $ is
    periodic.
\end{theorem}

\begin{cor}
Let $M$ be a graded finitely generated (finitely presented)
right module over a right Noetherian (respectively, right
coherent) connected graded algebra $A$. If $\gkd M \le 1$,
then the Hilbert function of $M$ is periodic.
\end{cor}

Notice that, in non-commutative projective geometry,
critical modules of linear growth are in one-to-one
correspondence of closed points to
 $\qcoh A$. So, a period of the Hilbert function becomes
 a numerical invariant of such point.

There exist finitely
generated algebras with non-periodic but bounded Hilbert
function, e.g., an algebra $k\la x,y | yxy,x^2,  x y^{2^n}
x, n \ge 1 \ra$; if the field $k$ is finite, this algebra
is even effective for series.

We do not know if any finitely presented algebra of linear
growth has periodic Hilbert function. However, we establish
the periodicity in several important cases.

\begin{cor}[Corollary~\ref{c-alg-GK1}]
\label{c-alg-GK1-intro}
 Let $A$  be a finitely generated connected algebra of Gelfand--Kirillov dimension one.
 Suppose that $A$ satisfies at least one of the following properties:

(i) two-sided or right Noetherian;

(ii) (semi)prime;

(iii) coherent;

(iv) finitely presented over a finite field;

(v) Koszul;

(vi) $A$ has finite Bakelin's rate.

Then the Hilbert function of $A$ is periodic.
\end{cor}

Here an algebra $A$ is said to be of finite Backelin's rate if
there is a number $r$ such that every space $\tor_i^A(k,k)$
is concentrated in degrees at most $r i$~\cite{brate}. $A$
is said to be (right graded) coherent if every finitely
generated homogeneous right--sided ideal is finitely
presented, or, equivalently, a kernel of any homogeneous
map $F_1 \to F_2$ of two finitely generated free modules
$F_1, F_2$ is finitely generated~\cite{bur, faith}.

An ideal in non-coherent algebras may also be finitely
presented and, moreover, it  may admit a free resolution of
finite type; to describe some of such ideals, we introduce
(in section~\ref{s-qcoh}) the following concept.

\begin{defi}
\label{def-qc-intro}
 Let $A$ be a connected graded algebra, and let
${\bf F}$ be a set of finitely generated  homogeneous
right-sided ideals in $A$. ${\bf F}$ is said to be
quasi-coherent family of ideals if $0 \in {\bf F}$ and  for
every $0 \ne I \in {\bf F}$ there are $J_1, J_2 \in {\bf
F}$ such that $J_1 \ne I, m_0(J_1) \le m_0(I)$, and $I/J_1
\cong A/J_2[-t]$ for some $t \in \Zpl$.

A quasi-coherent family ${\bf F}$ is said to be of degree
$d$ if $m_0(I) \le d$ for all $I \in {\bf F}$.
\end{defi}

If the maximal ideal $\overline A := A_1 \oplus A_2 \oplus
\dots \in {\bf F}$, a quasi-coherent family ${\bf F}$ is
called {\it coherent}~\cite{prepr}; for example, all
finitely presented monomial algebras admits coherent
families of finite degree~\cite{prepr}, as well as some
homogeneous coordinate rings~\cite{cnr}. A coherent family
of degree 1 is called {\it Koszul filtration}; Koszul
filtrations  have been studied in a number of
papers~\cite{i-kos, conca, con2, crv, ctv, kosf}. A
quasi-coherent family of degree one is called {\it Koszul
family}: such families exist, e.~g., in  homogeneous
coordinate rings of some finite sets of points in
projective spaces~\cite{pol}.

It is not hard to see that every ideal in a quasi-coherent
family admits free resolution of finite type. Using the
above results on the sets of Hilbert series, we deduce

\begin{prop}[Corollary~\ref{c-cq-fin}]
\label{p-cq-fin-intro} Let ${\bf F}$ be a quasi-coherent
family of degree $d$ in a finitely presented algebra $A$.
Then the set of Hilbert series of ideals $I \in {\bf F}$ is
finite.
\end{prop}

A version of Proposition~\ref{p-cq-fin-intro} for coherent
families has been proved in~\cite{prepr}; here we just
generalize it to quasi-coherent families using
Corollary~\ref{c-modules_best-intro}.

It is proved in~\cite{prepr, kosf} that every ideal in a
coherent family of finite degree has rational Hilbert
series, and that every ideal in a Koszul filtration has
also rational Poincare series.
Proposition~\ref{p-cq-fin-intro} allows us to establish
similar properties in quasi-coherent case.

\begin{cor}[Corollaries~\ref{c-cq-ratio},~\ref{c-koszfam_poincare}]
\label{c-ratio-intro} Let ${\bf F}$ be a quasi-coherent
family of degree $d$ in a finitely presented algebra $A$.

(i) For every ideal  $I \in {\bf F}$ there are two
polynomials $p(z), q(z)$ with integer coefficients such
that $I(z) = A(z)\frac{\displaystyle p(z)}{\displaystyle
q(z)}$.

(ii) Assume that $d=1$, i.~e. ${\bf F}$ is a Koszul family.
Then for every ideal  $I \in {\bf F}$ its Poincare series
$P_I(z) : = \sum_{i \ge 0} (\dim \tor_i(I,k)) z^i$ is a
rational function.
\end{cor}

\subsection{Acknowledgement}

I am grateful to Leonid Positselski, J. Tobias Stafford,
and Viktor Ufnarovski for helpful remarks and discussions.

\subsection{Notations and assumptions}

\label{ss_notations}

We will deal with ${\bf Z}_+$--graded connected associative
algebras over a fixed field $k$, that is, algebras of the
form   $R = \bigoplus_{i \ge 0} R_i$  with $R_0 = k$. By
assumption, all our modules and ideals are graded and
right-sided.

  For an $R$--module $M$, we will denote by $H_i M$
the graded vector space $\tor_i^R (M,k)$. By $H_i R$ we
will denote the graded vector space $\tor_i^R (k,k)  = H_i
k_A$. In particular, the vector space $H_1 R$ is isomorphic
to the $k$--span of a minimal set of homogeneous generators
of $R$, and  $H_2 R$ is isomorphic to the $k$--span of a
minimal set of its homogeneous relations. Analogously, the
space $H_0 M$ is the span of generators of $M$, and $H_1 M$
is the span of its relations.

Let $m(M) = m_0(M)$ denote the supremum of degrees of
minimal homogeneous generators of $M$: if $M$ is just a
vector space with trivial module structure, it is simply
the supremum of degrees of elements of $M$. For $i\ge 0$,
let us also put $m_i (M) := m (H_i M) = \sup \{ j |
\tor_i^R(M,k)_j \ne 0 \}$. Similarly, let us put $m_i(R) =
m (H_i R) = m_i (k_R)$. For example, $m(R) = m_0(R)$ is the
supremum of degrees of the generators of $R$, and $m_1 (R)$
(respectively, $m_1 (M)$) is the supremum of degrees of the
relations of $R$ (resp., of~$M$).


Note that the symbols $H_i R$ and  $m_i R$ for an algebra
have different meaning that  the respective symbols $H_i
R_R$ and $m_i R_R$ for $R$ considered as a module over
itself; however, the homologies $H_i R_R$ are trivial, so
that there is no place for confusion.


For a graded locally finite vector space (algebra,
module...) $V$, its Hilbert series is defined as the formal
power series $V(z) = \sum_{i \in {\bf Z}}  (\dim V_i) z^i$.
For example, the Euler characteristics of a minimal free
resolution of the trivial module $k_R$ leads to the formula
\begin{equation}
\label{euler}
    R(z)^{-1} = \sum_{i \ge 0} (-1)^i H_i  R(z). \end{equation}
As usual, we write $\sum _{i \ge 0} a_i z^i = o(z^n)$ iff
$a_i =0$ for $i \le n$.

Let us introduce a lexicographical total order on the set
of all power series with integer coefficients, i.e., we put
$\sum_{i \ge 0} a_i z^i >_{lex} \sum_{i \ge 0} b_i z^i$ iff
there is $q \ge 0$  such that $a_i = b_i$  for $i < q$ and
$a_q > b_q$. This order extends the coefficient-wise
partial order  given by $\sum_{i \ge 0} a_i z^i \ge \sum_{i
\ge 0} b_i z^i$ iff $a_i \ge b_i$ for all $i \ge 0$.

\section{Properties of sets of Hilbert series}

\label{s_sets_Hilbert}

 The following theorem of Anick shows that the set of
  Hilbert series of algebras of bounded (by degrees and numbers)
  generators and relations is well-ordered.

\begin{theorema}[{\cite[Theorem~4.3]{an4}}]
\label{ant} Given three integers $n,a,b$, let $C(n,a,b)$ be
the set of all $n$--generated connected algebras $R$ with
$m_1(R) \le a$ and $m_2(R) \le b$ and let ${\mathcal H}
(n,a,b)$ be the set of Hilbert series of such algebras.
Then the ordered set $({\mathcal H} (n,a,b) ,  >_{lex} )$
admits no infinite ascending chains.
\end{theorema}

The example of an infinite {\it descending} chain of
Hilbert series in the set $C(7,1,2)$ is constructed
in~\cite[Example~7.7]{an4}.

We will prove this theorem in a more general form, with
(bi)modules instead of algebras. The proof is bases on the
same idea as the original proof in~\cite{an4}.

\begin{theorema}
\label{an-mod} \label{mod-gen_a}
Let $D \ge D_0$ be two integers. Consider a set $E(D_0, D)$
of all graded (bi)modules $M$ over connected graded
algebras $A$ such that  $A$ and $M$ are generated by $\le
D$ elements, the generators of $M$ have degrees at least
$D_0$, and the generators and relations of both $A$ and $M$
have degrees at most $D$. Then the set of all Hilbert
series of (bi)modules from $E(D_0, D)$ has no infinite
ascending chains.
\end{theorema}

Let us first introduce an additional notation.

Given four formal power series $V(z), R(z), W(z), S(z)$,
let us denote by  ${\mathcal M} = {\mathcal M}(V(z), R(z),
W(z), S(z))$
 the set of all modules $M$ over algebras $A$ such that
 $H_1 A(z) = V(z)$, $H_2 A(z) = R(z)$, $H_0 M(z) = W(z)$,
 and $H_1 M(z) = S(z)$.

We may assume that all these algebras are generated by the
same vector space $V$ and all our
modules are generated by the  same vector space $W$. Let $
D = \max{ \deg R(z), \deg S(z)}$, and let $m= \dim
T(V)_{\le D} $, $n = \dim (W \otimes T(V))_{\le D}$. Since
every relation $r \in H_2 A$ of an algebra $A$ is an
element of $T(V)_{\le D} $, it may be considered as an
element of $k^m$; so, the vector space $H_2 A$ is uniquely
determined by a vector in $k^{rm}$, where $r = R(1)$ is the
dimension of $H_2 A$. Analogously, every vector space of
relations of a module $M \in {\mathcal M}$  is uniquely
determined by a vector  in $k^{sn}$, where $s = S(1)$ is
the dimension of $H_1 M$. This means that every module $M =
M_u \in {\mathcal M}$ is uniquely determined by a vector
$u$ in a vector space $Q = k^{rm+sn}$.

Consider a topological space $\overline{\mathcal M}(V(z),
R(z), W(z), S(z)) := \{ u | M_u \in {\mathcal M} \} \subset
Q$ with induced Zarisski topology.

\begin{lemma}
\label{lem_zarisski_closed}
 Let $V(z), R(z), W(z), S(z)$ be four polynomials
with positive integer coefficient, and let $h(z)$ be a
formal power series. Then two subsets  $L_{>}(h(z)) := \{ u
| M_u(z) \ge h(z) \}$ and $L_{>_{lex}}(h(z)) :=   \{ u |
M_u(z) {\ge}_{lex} h(z) \}$ in $\overline{\mathcal M}(V(z),
R(z), W(z), S(z))$ are closed and algebaric.
\end{lemma}

\begin{proof}[Proof of Lemma~\ref{lem_zarisski_closed}]
Let $M = M_u \in {\mathcal M}$ be a
module
over an algebra $A$, and
let $R,S$ be minimal sets of relations of $A$ and $M$. Then
$M$ is a quotient of the free
module
$F = W \otimes T(V) $
 by a
$T(V)$-submodule
$N = W \otimes T(V) R T(V) + S T(V)$. Put $h(z) =
\sum_{i\ge 0} h_i z^i$ and $\tilde h(z) = \sum_{i \ge 0}
h_i z^i := F(z) - h(z) $.

The condition $u \in L_{>}(h(z))$ means that $N(z) \le
\tilde h(z) $, that is, $\dim N_i \le h_i$  for every $i
\ge 0$. For every $i \ge 0$, the last condition means that
the rank of the vectors generating the vector space $N_i$
is bounded above by $h_i$. Obviously, this condition is
algebraic for $u$, because it simply means suitable minor
determinants vanish. Therefore,  the set $L_{>}(h(z))$ is a
countable intersection of closed subsets, hence it is
closed.

Now, the condition $u \in L_{>_{lex}}(h(z))$ means that
$N(z) \le_{lex} \tilde h(z)$. Hence, the set
$L_{>_{lex}}(h(z))$ is a countable intersection of the sets
$L_i$, $i \ge 0$, where $L_i = \{ u | N(z) \le h(z) +
o(z^i) \} \bigcup_{j=0}^i \{ u | \dim N_{j} < h_j \}
\subset \overline{\mathcal M}(V(z), R(z), W(z), S(z)) $.
Every set $L_i$ is algebraic, hence $L$ is algebraic as
well.
\end{proof}

\begin{proof}[Proof of Theorem~\ref{an-mod}]
First, up to a shift of grading we may assume that $D_0 =
0$. Second, every bimodule over an algebra $A$ may also be
considered as a right module over the algebra $B = A
\otimes A^{op}$. Since $m_1(B) = m_1(A) \le D$ and $m_2(B)
\le \max \{ m_1(A)^2 , m_2 (A) \} \le D^2$, it is
sufficient to prove Theorem~\ref{an-mod} (for every $D$)
for a subset $E'(0, D) \subset E(0,D) $ which consists of
right modules (that is, of bimodules with zero left
multiplication). Indeed, because the relations of a right
module $M$ as a bimodule may have degrees at most $\max \{
m_0(M)+m_1 (A) , m_1 (M) \} \le 2D$, the statement for the
set $E(0,D)$ will follow from the same statement for the
set $E'(0, \max \{ 2D, D^2 \} )$.


Assume that there is an infinite ascending chain
$$
   M_1(z) <_{lex}    M_2(z) <_{lex} \dots
$$
in $E'(0, D)$ of modules over some algebras $A_1, A_1,
\dots \in C(D,D,D)$. We may assume that all these algebras
are generated by the same finite-dimensional  graded vector
space $V$ and that their minimal vector spaces of relations
$R_1, R_2, \dots \subset T(V) $ have the same Hilbert
series (polynomial) $R_i(z) = R(z)$. Analogously, we assume
that all modules $M_i$ are generated by the same
finite-dimensional  graded vector space $W$ and has the
minimal vector spaces of relations $S_1, S_2, \dots $ of
the same finite Hilbert series
 $S_i(z) = S(z)$.

Thus we obtain an infinite descending set of closed subsets
in $\overline{\mathcal M}(V(z), R(z), W(z), S(z))$:
$$
L_{>_{lex}}(M_1(z)) \supset L_{>_{lex}}(M_2(z)) \supset
\dots,
$$
a contradiction.
\end{proof}

The following theorem is the main result of this section.
Another its proof may be found in~\cite{prepr}.

\begin{theorem}

\label{main}

 Let $n,a,b,c, m, p_1, p_2, q, r$  be 9 integers.

(a)Let $D(n,a,b,c)$ denote the set of all connected
algebras $A$ over a fixed field $k$ with at most $n$
generators such that $m_1 (A) \le a, m_2(A) \le b$, and
$m_3(A) \le c$. Then the set of Hilbert series of algebras
from $D(n,a,b,c)$ is finite.

(b)   Let $DM =DM(n,a,b,c, m, p_1, p_2, q, r)$ denote the
set of all graded right modules  over algebras from
$D(n,a,b,c)$ with at most $m$ generators such that  $M_i =
0$ for $i<p_1$, $m_0(M) \le p_2, m_1(M) \le q$, and $m_2(M)
\le r$. Then the set $HDM$ of Hilbert series of modules
from $DM$ is finite.
\end{theorem}

For this statement, we need the following standard version
of Koenig lemma.
\begin{lemma}
\label{Koenig} Let $P$ be a totally ordered set satisfying
both ACC and DCC. Then $P$ is finite.
\end{lemma}

\begin{proof}[Proof of Theorem~\ref{main}]
 Up to a shift of grading, we may assume
that $p_1 = 0$. Let $M \in DM$ is a module over an algebra
$A \in D(n,a,b,c)$. Consider the first terms of the minimal
free resolution of $M$:
$$
 0 \to \Omega \to H_0 (M) \otimes A  \to M \to 0.
$$
Here the syzygy module $\Omega$ has generators in degrees
at most $m_1(M) \le q$ and relations in degrees at most
$m_2 (M) \le r$. Since $H_0 (\Omega) = H_1(M)$, the number
$\dim H_0 (\Omega)$ of its generators is not greater than $
\dim (H_0 (M) \otimes A)_{\le q} \le m (1+ n +\dots + n^q)
=:D'$. We have that $\Omega \in E(0, D)$ for $D = \max \{
D', q,r \}$, so, the set of all possible Hilbert series for
$\Omega$ satisfies ACC.

Assume for a moment that $M=A_1 \oplus A_2 \oplus \dots$
(to connect the cases~$(a)$ and~$(b)$, we take here $m=n,
p_2 = a, q=b$, and $r=c$). We have here $M(z) = A(z)-1$ and
$H_0(M) = H_1 (A)$. Taking Eulerian characteristics, we
obtain
 $\Omega (z) - H_1(A)(z) A(z) + M(z) =0$, or
$$
       A(z) = (1-\Omega (z))(1-H_1(A)(z))^{-1} = (1-\Omega (z))(1+H_1(A)(z)+H_1(A)(z)^2 + \dots).
$$
Notice that the order $ <_{lex} $ is compatible with the
multiplication by a formal power series with positive
coefficients. Because there is only finite number of
possibilities for $H_1(A)(z)$, the set
of the
Hilbert series $A(z)$ of that type satisfies DCC. By
Theorem~\ref{ant}, this set also  satisfies ACC. In the
view of Lemma~\ref{Koenig}, the statement~$(a)$ follows.

Now, return to the general case~$(b)$. We have $M(z) = H_0
(M)(z) A(z) - \Omega (z)$. Here there is only finite number
of possibilities for $H_0 (M)(z)$ because of dimension
arguments; also, there is only finite number of
possibilities for $A(z)$ by part~$(a)$. It follows that the
set $HDM$ of all such Hilbert series $M(z)$ satisfies DCC.
Because $HDM \subset E(0, \max \{ n,a,b,m,p_2,q\})$, it
also satisfies ACC by Theorem~\ref{an-mod}. By
Lemma~\ref{Koenig}, it if finite.
\end{proof}

If we restrict our consideration to the modules over a
single algebra $A$, the additional condition $m_3(A) <
\infty$ may be omitted.

\begin{cor}
\label{modules_best} Let $A$ be a finitely presented
algebra. Given five integers $m, p_1, p_2, q, r$, consider
a set $DM_A = DM_A(m, p_1, p_2, q, r)$ of finitely
presented $A$-modules $M$ with at most $m$ generators such
that  $M_i = 0$ for $i<p_1$, $m_0(M) \le p_2, m_1(M) \le
q$, and $m_2(M) \le r$. Then the set of Hilbert series of
modules from $DM_A$ is finite.
\end{cor}

\begin{proof}
Like the proof of Theorem~\ref{main}, we  may assume that
$p_1 =0$ and consider the exact sequence
$$
 0 \to \Omega \to H_0 (M) \otimes A  \to M \to 0.
$$
As before, we see that the set of all possible Hilbert
series for $\Omega$ satisfies ACC, so, the set $HDM_A$ of
Hilbert series $M(z) = H_0 (M)(z) A(z) - \Omega (z)$
satisfies DCC. Because $DM_A \subset E(0, D)$ for $D =
\max\{ m_1(A), m_2(A), \dim H_0 (A), p_2, q_r  \}$,
Lemma~\ref{Koenig} implies that the set of Hilbert series
of modules from $DM_A$ is finite.
\end{proof}

\begin{cor}
\label{mod} Let $D > 0$ be an integer, and let $A$ be a
finitely presented algebra.

(a) The set of Hilbert series of (two-sided or right-sided)
ideals in $A$ generated in degrees at most $D$ satisfies
DCC.

(b) The set of Hilbert series of right-sided ideals in $A$
having generators and relations in degrees at most $D$ is
finite.

\end{cor}

\begin{proof}
Just apply Theorem~\ref{an-mod} and
Corollary~\ref{modules_best} to the modules $A/I$, where
$I$ is an ideal.
\end{proof}

\section{Periodic Hilbert functions}

\label{s-periodic}

\begin{theorem}
\label{main-GK1} \label{t-main-GK1}
 Let $M$ be a finitely presented graded module over a connected finitely presented algebra $A$.
Suppose that $\gkd M =1$ and  at least one of the following
conditions holds:

  (a)  the field $k$ is finite;

  (b) the vector spaces
  $\tor_2^A (M,k)$  and $\tor_3^A (k,k)$ have finite dimensions.

    Then the Hilbert series $M(z)$ is rational, that is, the Hilbert function $f_M(n) = \dim M_n $ is
    periodic.
\end{theorem}

\begin{cor}
\label{c-alg-GK1}
 Let $A$  be a finitely generated algebra of Gelfand--Kirillov dimension one.
 Suppose that $A$ satisfies at least one of the following properties:

(i) two-sided or right Noetherian;

(ii) (semi)prime;

(iii) coherent;

(iv) finitely presented over a finite field;

(v) Koszul;

(vi) of finite Bakelin's rate.

Then the Hilbert funcion of $A$ is periodic.
\end{cor}

\begin{proof}[Proof of Corollary~\ref{c-alg-GK1}]
Recall that any affine algebra $A$ with $\gkd A = 1$ is
PI~\cite{ssw}, hence it has rational Hilbert series
provided that it is Noetherian~\cite{Lorenz}. Moreover, it
is shown in~\cite{Lorenz} that every Noetherian module over
a PI algebra has rational Hilbert series; since every
Noetherian $A$--bimodule is a Noetherian module over the
algebra $A^{op} \otimes A$, it follows that any weak
Noetherian (i.~e., satisfying ACC for two-sided ideals)
 PI algebra has rational Hilbert series.
 This proves the case~$(i)$.

Any semiprime algebra of Gelfand--Kirillov dimension one is
a finite module over its Noetherian center~\cite{ssw},
hence it has periodic Hilbert function. This proves~$(ii)$.

The rest 5 cases  follow from Theorem~\ref{t-main-GK1}.
Notice that the case~$(v)$ has been also proved by L.
Positselski (unpublished).
\end{proof}

\begin{lemma}
\label{lem1} Let  $M$  be an infinite--dimensional graded
(by nonnegative integers) module over a connected graded
algebra $A$. Let $M^n$ denote the right $A$--module
$\bigoplus\limits_{t \ge n} M^t$ with degree shifted by
$n$, i.e., $(M^n)_t = M_{n+t}$ for all $t \ge 0$. Then the
following conditions are equivalent:

(a)  $\gkd M = 1$ and $M(z)$ is a rational function;

(b) for some $i \ne j$, we have  $M^i (z) = M^j(z)$; 

(c)  the set of Hilbert series ${\mathcal H}_M = \{ M^n (z)
| n \ge 0\}$  is finite.


In this case, the sequence $\{ \dim M^n   \}$ is periodic
with a period $d$ such that $d \le |{\mathcal H}_M|$ and $d
\le |i-j|$.
\end{lemma}

\begin{proof}
If $\gkd M = 1$, the rationality of $M(z)$ means that there
are two positive integers $d, D$ such that $\dim M_n = \dim
M_{n+d}$ for all $n \ge D$, that is, $M^n(z) = M^{n+d}(z)$
for all $n \ge D$. In this case the set ${\mathcal H}_M =
\{ M^n (z) | n < D+d \}$ is finite.

On the other hand, if   $M^n(z) = M^{n+d}(z)$ for some $n
\ge 0, d > 0$, then $\dim M_i = \dim M_{i+d}$ for all $i
\ge n$, so, $\dim M_i \le \max_{j \le n+d} M_j$. In
particular, it follows that  $\gkd M = 1$.

Finally, if the set ${\mathcal H}_M$ is finite, then $(b)$
obviously holds.
\end{proof}

\begin{rema}
By the same way, it may be shown more: if $\gkd M \le 1$
and $M^i (z) \le M^j(z)$ for some  $i \ne j$, then the
Hilbert function of $M$ is periodic.
\end{rema}

\begin{lemma}
\label{Tor_M_i} Let  $M$  be a nonegatively graded  module
over a connected graded algebra $A$, and let $M^n$ be as in
Lemma~\ref{lem1}. Then $m_i(M^n) < \max\{ m_{i}(M),
m_{i+1}(A)$ for all $n \ge 1, i \ge 0$.
\end{lemma}

\begin{proof}
By definition, $m_i(M^0) =  m_{i}(M)$ for all $i \ge 0$.
Let us prove by induction on $n \ge 1$ that $\tor_{i}^A
(M^{n+1},k)_{j-1} = 0 $ provided that $\tor_{i}^A
(M^{n},k)_j = \tor_{i+1}^A (k,k)_j =0$.

Let $s = \dim A_{n}$. The exact triple
$$
    0 \to M^{n+1}[-1] \to M^{n} \to k^{s} \to 0
$$
leads to the triple
$$
     \tor_{i+1}^A (k^{s} ,k) \to  \tor_{i}^A (M^{n+1},k) [-1] \to  \tor_{i}^A (M^{n},k)
$$
for every $i \ge 0$.   It remains to notice that
$\tor_{i+1}^A (k^{s} ,k) = \bigoplus_{1}^s  \tor_{i+1}^A
(k,k)$.
\end{proof}

\begin{proof}[Proof of Theorem~\ref{main-GK1}]
Up to the shift of grading, we may assume that $M_i = 0$
for $i<0$. Put $g_A = m(A), g_M = m(M), r_A = m_2(A)$, and
$r_M = m_1(M)$. Let $N$ be a number such that $\dim M_n \le
N$ for all $n \ge 0$. By Lemma~\ref{Tor_M_i}, every module
 $M^n$  is isomorphic to a quotient of a free module
 $F = V \otimes A$  by a submodule generated by a
homogeneous subspace $W \subset  F_{[1..r_M]} = F_1 \oplus
\dots \oplus F_{r_M}$, where $V = V_0 \oplus \dots \oplus
V_{g_M}$ is a graded vector space with $\dim V_i = N, 0 \le
i \le g_M$.

{\it Case $(a)$.} Let $q$ be the cardinality of $k$. Since
$\dim F_{[1..r_M]}$  is finite (namely, $\dim F_{[1..r_M]}
< g_A r_M N =: T$, so $| F_{[1..r_M]}| < q^{T} $), then
there is only finite number of possibilities for its subset
$W$ (namely, there are at most $2^{q^{T} }$ proper subsets
of $F_{[1..r_M]}$). Then there is a finite number of the
isomorphism types of the modules $M^n$, and so, by
Lemma~\ref{lem1}, the sequence $\{ \dim M_n   \}$ is
periodic with a period $d < 2^{q^{T} } $.

{\it Case $(b)$.} By Lemma~\ref{Tor_M_i}, for every $n \ge
2$ we have $m_2(M^n) \le m_3(A) -1$. By
Corollary~\ref{modules_best}, there is only finite number
of possibilities for Hilbert series $L(z)$ of modules $L$
with bounded $\dim H_0 (L)$ and $m_i(L)$ for $i \le 3$.
Thus we conclude that the set of all Hilbert series of the
modules $M^i$ is finite. It remains to apply
Lemma~\ref{lem1}.
\end{proof}

\section{Quasi-coherent families of ideals}

\label{s-qcoh}

\begin{defi}
\label{def-qc}
 Let $A$ be a connected graded algebra, and let
${\bf F}$ be a set of finitely generated  homogeneous
right-sided ideals in $A$. ${\bf F}$ is said to be
quasi-coherent family of ideals if $0 \in {\bf F}$ and for
every $0 \ne I \in {\bf F}$ there are $J_1, J_2 \in {\bf
F}$ such that $J_1 \ne I, m_0(J_1) \le m_0(I)$, and $I/J_1
\cong A/J_2[-t]$ for some $t \in \Zpl$.

A quasi-coherent family ${\bf F}$ is called of degree $d$
if $m_0(I) \le d$ for all $I \in {\bf F}$.
\end{defi}

A quasi-coherent family of degree 1 is called {\it Koszul
family} of ideals. It was introduced by
Polishchuk~\cite{pol} in order to prove that homogenous
coordinate rings of some sets of points in projective
spaces are Koszul. If $A_{\ge 1} \in {\bf F}$, then a
quasi-coherent family is called {\it
coherent}~\cite{prepr}. The term "quasi-coherent family"
appears because all ideals in such a family are finitely
presented, like finitely generated submodules  in a
quasi-coherent module. Moreover, there is

\begin{prop}
\label{p-qcoh} Let ${\bf F}$ be a quasi-coherent family in
an  algebra $A$. Then every ideal $I \in {\bf F}$ has free
resolution of finite type.  If ${\bf F}$ has degree $d$,
then for every $I \in {\bf F}$ we have $m_i(I) \le m_0(I)
+id $ for all $i \ge 0$.
\end{prop}

The proof is the same as for coherent families
in~\cite{prepr}.

\begin{proof}
We proceed by induction in $i$ and in $I$ (by inclusion of
ideals $J_1$ with $m_0(J_1) \le m_0(I)$). Let $J_1, J_2, t$
be as in Definition~\ref{def-qc}; in particular, $t \le
m_0(I) \le d$. The exact sequence
$$
           0 \to J_1 \to I \to A/J_2[-t] \to 0
$$
leads for every $i \ge 1$ to the following fragment of the
long exact sequence of $\tor$'s:
$$
             \dots \to H_i (J_1)_j \to H_i (I)_j \to
             H_{i-1}(J_2)_{j-t} \to \dots
$$
By induction, we have $m_{i}(I) \le \max \{ m_i(J_1) ,
m_{i-1} (J_2) +t \} < \infty$. If $t \le d$, we have  also
$m_{i}(I) \le t + id$.
\end{proof}

\begin{cor}
\label{c-cq-fin} Let ${\bf F}$ be a quasi-coherent family
of degree $d$ in a finitely presented algebra $A$. Then the
set of Hilbert series of ideals $I \in {\bf F}$ is finite.
\end{cor}

\begin{proof}
In the view of Proposition~\ref{p-qcoh}, we can apply
Corollary~\ref{mod}.
\end{proof}

\begin{cor}
\label{c-cq-ratio} Let ${\bf F}$ be a quasi-coherent family
of degree $d$ in a finitely presented algebra $A$. Then for
every ideal  $I \in {\bf F}$ there are two polynomials
$p(z), q(z)$ with integer coefficients such that $I(z) =
A(z)\frac{\displaystyle p(z)}{\displaystyle q(z)}$.
\end{cor}

\begin{proof}
Because of Corollary~\ref{c-cq-fin}, we can apply exactly
the same arguments as for coherent families,
see~\cite[Theorem~4.5]{prepr}.
\end{proof}

\begin{cor}
\label{c-koszfam_poincare} Let ${\bf F}$ be a Koszul family
(that is, a quasi-coherent family of degree $1$) in a
finitely presented algebra $A$. Then for every ideal  $I
\in {\bf F}$ its Poincare series $P_I(z) : = \sum_{i \ge 0}
(\dim \tor_i(I,k)) z^i$  is a rational function.
\end{cor}

\begin{proof}
Every ideal in a Koszul family is a Koszul module
(see~\cite{pol} or Proposition~\ref{p-qcoh}). By
Corollary~\ref{c-cq-ratio}, we have $I(z) = p(z) / q(z)$.
By Koszulity, we have $I(z) = P_I(-z) A(z)$,  so that $P(z)
= p(-z)/q(-z)$.
\end{proof}

\end{document}